\documentclass[12pt]{amsart}

\usepackage[utf8]{inputenc}
\usepackage{amsmath}
\usepackage{graphics}
\usepackage{amsthm, multicol}
\usepackage{bm}
\usepackage{filecontents}
\usepackage{amssymb}
\usepackage{setspace, tikz, float,graphicx, tikz-cd}
\usepackage{mathrsfs}
\usepackage{mathabx}
\usepackage{interval}
\usepackage{amsfonts,latexsym,amscd,euscript,graphicx}
\providecommand{\keywords}[1]{\textbf{\textit{Index terms---}} #1}
\usepackage{framed, mathtools}
\usepackage{fullpage}
\usepackage{hyperref}
\usepackage[margin=0.8 in]{geometry}
\usepackage{setspace}
\usepackage{enumerate}

\newcommand{\N}{{\mathbb N}}

   \hypersetup{colorlinks=true,citecolor=blue,urlcolor =black,linkbordercolor={1 0 0}}

\allowdisplaybreaks[1]
\usepackage{thmtools}
\declaretheoremstyle[notefont=\bfseries,notebraces={}{},%
    headpunct={},postheadspace=1em,spaceabove=0.5em,spacebelow=0.5em]{mystyle}
\declaretheorem[style=mystyle,numbered=no,name=Theorem]{thm-hand}
\declaretheorem[style=mystyle,numbered=no,name=Conjecture]{conj-hand}
\declaretheorem[style=mystyle,numbered=no,name=Definition.]{defn-hand}
\declaretheorem[style=mystyle,numbered=no,name=Theorem.]{thm-no}
\declaretheorem[style=mystyle,numbered=no,name=Conjecture.]{conj-no}
\declaretheorem[style=mystyle,numbered=no,name=Lemma.]{lemma-no}
\declaretheorem[style=mystyle,numbered=no,name=Lemma]{lemma-hand}
\makeatletter
\def\resetMathstrut@{%
  \setbox\z@\hbox{%
    \mathchardef\@tempa\mathcode`\[\relax
    \def\@tempb##1"##2##3{\the\textfont"##3\char"}%
    \expandafter\@tempb\meaning\@tempa \relax
  }%
  \ht\Mathstrutbox@\ht\z@ \dp\Mathstrutbox@\dp\z@}
\makeatother
\begingroup
  \catcode`(\active \xdef({\left\string(}
  \catcode`)\active \xdef){\right\string)}
  \endgroup
\mathcode`(="8000 \mathcode`)="8000

\theoremstyle{definition}

\usepackage{listings}
\theoremstyle{remark}

\newcommand{\e}{\epsilon}

\newcommand{\lb}{\langle}
\newcommand{\rb}{\rangle}

\newcommand{\nc}{\newcommand}

\nc{\on}{\operatorname}
\nc{\Spec}{\on{Spec}}
\vspace{-15mm}
\title{On Symmetric But Not Cyclotomic Numerical Semigroups} 

\author{Mehtaab Sawhney}
\thanks{Massachusetts Institute of Technology, Cambridge MA. Email: \texttt{msawhney98@yahoo.com}}
\author{David Stoner}
\thanks{Harvard University, Cambridge MA. Email: \texttt{dstoner@college.harvard.ed}}
\date{\today}

\begin{document}

\maketitle
\vspace{-10mm}
\begin{abstract}
A numerical semigroup is called \textit{cyclotomic} if its corresponding numerical semigroup polynomial $P_S(x)=(1-x)\sum_{s\in S}x^s$ is expressable as the product of cyclotomic polynomials. Ciolan, Garc{\'\i}a-S{\'a}nchez, and Moree conjectured that for every embedding dimension at least $4$, there exists some numerical semigroup which is symmetric but not cyclotomic. We affirm this conjecture by giving an infinite class of numerical semigroup families $S_{n, t}$, which for every fixed $t$ is symmetric but not cyclotomic when $n\ge \max(8(t+1)^3,40(t+2))$ and then verify through a finite case check that the numerical semigroup families $S_{n, 0}$, and $S_{n, 1}$ yield acyclotomic numerical semigroups for every embedding dimension at least $4$. 
\end{abstract}

\section{Background}
A subset $S\subseteq \mathbb{Z}_{\ge 0}$ is a numerical semigroup if $S$ contains zero, is closed under addition, and contains all but finitely many nonnegative integers. Every numerical semigroup $S$ has a unique minimal generating set, written as $\lb n_1, n_2, \ldots, n_e\rb$ and we will identify $S=\lb n_1, n_2, \ldots, n_e\rb$. We say that $e$ is the \textit{embedding dimension} of $n$. Given a numerical semigroup $S$, we define the \textit{Hilbert series} $H_S(x)$ to be the formal power series $H_S(x)=\sum_{s\in S}x^s$. Furthermore we define the \textit{numerical semigroup polynomial} $P_S(x)$ by $P_S(x)=(1-x)H_S(x)$. This is a bona fide polynomial as $S$ has finite complement.

We are interested in two particular properties of these numerical semigroup polynomials. A numerical semigroup $S$ is \textit{symmetric} if $P_S(x)$ has palindromic coefficients; that is, $P_S(x)=x^kP_S(x^{-1})$ for some nonnegative integer $k$. For example, $P_{\lb 5, 6, 7, 8\rb}(x)=x^{10}-x^9+x^5-x+1$, so $\lb 5, 6, 7, 8\rb$ is symmetric. A numerical semigroup $S$ is called \textit{cyclotomic} if every root of $P_S(x)$ lies on the unit circle. Due to a result of Kronecker's (see, e.g. \cite{Kronecker}) this condition is equivalent to the requirement that $P_S(x)$ is the product of cyclotomic polynomials. For each positive integer $e$, let $SY_e$ and $CY_e$ denote the set of symmetric and cyclotomic numerical semigroups, respectively, which have embedding dimension $e$. An introduction to cyclotomic polynomials and numerical semigroups is given by Ciolan, Garc{\'\i}a-S{\'a}nchez, and Moree in \cite{Moree}. See also \cite{1} and \cite{2} for further discussion on numerical semigroup polynomials. 

Since all cyclotomic polynomials are symmetric, it follows that for every positive integer $e$ we have $CY_e\subseteq SY_e$. The following is also known, and noted in \cite{Moree}:
\begin{thm-no}\cite{Moree}
For $e\in \{1, 2, 3\}$, a numerical semigroup of embedding dimension $e$ is cyclotomic if and only if it is symmetric; that is, $CY_e=SY_e$. 
\end{thm-no}
Despite this equivalence for small dimensions, Ciolan, Carcia-S{\'a}nchez, and Moree \cite{Moree}, proposed the following conjecture:
\begin{conj-no}
For every positive integer $e\ge 4$, there exists a numerical semigroup of embedding dimension $e$ which is symmetric but not cyclotomic. That is, $CY_e\neq SY_e$. 
\end{conj-no}
Our main result strengthens this conjecture to a significant extent by giving an infinite family of obstructions. This conjecture was also proved independently by different methods by Herrera-Poyatos and Moree. \cite{moree2}
\begin{defn-hand}
For every nonnegative integer $t$ and positive integer $n\ge 6t+2$, let $S_{n, t}$ be the numerical semigroup generated by:
\begin{itemize}
\item The $t$ pairs $n-2t+4i, n-2t+4i+1$.  for $0\le i\le t-1$,
\item The $n-6t-1$ integers $n+2t, \ldots, 2n-4t-2$, and
\item The $t$ integers $2n-4t+4j-1$ for $0\le j\le t-1$.
\end{itemize}
\end{defn-hand}
\begin{thm-hand}[1.]
For every nonnegative integer $t$, and for every positive integer $n\ge \max$ ($8(t+1)^3,$ $40(t+2)$), the numerical semigroup $S_{n, t}$ has embedding dimension $n-3t-1$, and is symmetric but not cyclotomic.
\end{thm-hand}
Furthermore by explicitly checking remaining dimensions, the families $S_{n, 0}$, and $S_{n, 1}$ yield symmetric but not cyclotomic numerical semigroups for all embedding dimensions at least $4$. Note that Ciolan, Carcia-S{\'a}nchez, and Moree \cite{Moree} proposed that $\lb m , m + 1, q m + 2q + 2, \ldots , q m + (m -1)\rb$ for $q\ge 2m+4$ would give a family of numerical semigroups that are symmetric but not cyclotomic. The result given for $S_{n,1}$ corresponds precisely to the case $q=1$ in the above conjecture. Therefore our results imply the following:
\begin{thm-hand}[2.]
For every positive integer $e\ge 4$, there exists a pair of numerical semigroups of embedding dimension $e$ which is symmetric but not cyclotomic. Therefore $CY_e\neq SY_e$. 
\end{thm-hand}
This conjecture was proved independently by Herrera-Poyatos and Moree in \cite{moree2} by different methods. 
\\
\\ Note that our results demonstrate that in fact $\displaystyle\lim_{e\to\infty}|SY_e \setminus CY_e|=\infty$.

\section{Results}
We first proceed with the proof of Theorem 1. It is possible to explicitly compute the Hilbert series and numerical semigroup polynomial of the $S_n$.
\begin{lemma-hand}[3.]
For each positive integer $n\ge 6t+1$, we have: 
\begin{align*}
P_{S_{n,t}}(x)&=x^{2n}-x^{2n-1}+x^{n-2t}(\frac{x^{4t+2}+1}{x^2+1})-x+1.
\\&=x^{2n}-x^{2n-1}+\sum_{i=0}^{2t}(-1)^ix^{n+2t-2i}-x+1
\end{align*} In particular, $P_{S_{n,t}}(x)$ is symmetric. 
\end{lemma-hand}
\begin{proof}
This verification is equivalent to showing that, apart from the generators, the set $S$ contains exactly:
\begin{itemize}
\item The integers $2n-4t+4j, 2n-4t-4j+1, 2n-4t-4j+2$ for $0\le j\le t-1$.
\item All positive integers $\ge 2n$. 
\end{itemize}
Note that in the case $t=0$, the set of generators of $S_{n,0}$ is $n,n+1,\ldots, 2n-2$. and therefore it suffices to demonstrate $\N /S_{n,0}=\{1,2,\ldots, n-1, 2n-1\}$. This follows as all other nonnegative integers lie in the interval $[kn, k(2n-2)]$ for some nonnegative integer $k$. For $t\ge 1$ we first prove that $S$ contains these integers, and then that $S$ does not contain other integers. Indeed, every integer in the first group is expressible as $n-2t$ or $n-2t+1$ plus some generator from one of the $t$ pairs. For integers of the second type, we prove by induction the following:
\\
\\
\textbf{Claim:} For every $k\ge 2n$, there is at most one integer in the range $[n+2t, k]$ which is not in $S$. 
\begin{proof}
The proof is via induction. The base case follows from $2n=(n-2t)+(n+2t)$ and the already shown members of $S$, excluding $2n-1$. For the inductive step, suppose that the result holds for $k=\ell\ge 2n$. Then $\ell+1-(n-2t)$ and $\ell+1-(n-2t+1)$ are both in $[n+2t, \ell]$. In particular, one of them is in $S$, so $\ell+1$ is in $S$ as well. The induction is complete.  
\end{proof}
Finally, we check that the rest of the positive integers are not in $S$. That is, every element in $S$ that is not a generator is at least $n+2t-1$, and includes $2n-1$. The first follows from the bound $2(n-2t)>n+2t>n+2t-1$. For the second, we have $3(n-2t)>2n>2n-1$, so $2n-1$ could only potentially be expressed as the sum of two generators. For this to be the case, due to size restrictions, both generators must be in the $t$ pairs $n-2t+4i, n-2t+4i+1$.  for $0\le i\le t-1$. But the sum of any two generators of this form is equivalent to $2n, 2n+1,$ or $2n+2 \mod 4$. Hence $2n-1$ is not expressible, and the lemma is proved.
\end{proof}
To demonstrate that this polynomial is not the product of cyclotomic polynomials when $n\ge c_t$, we prove that $P_{S_{n,t}}(x)$ has at most $2n-1$ roots which are roots of unity (including multiplicity).

\begin{lemma-hand}[4.]
Consider a smooth function $g(x)$ on the interval $[a,b]$ with no triple roots. Then the number of roots of $g'(x)$ at least one less than the number of roots of $g(x)$ including multiplicity.
\end{lemma-hand}
\begin{proof}
Suppose that $g(x)$ has roots $r_1<\ldots<r_m$ with multiplicity at most $2$ for each $r_i$. It follows that $g'(x)$ has at least one root in the interval $(r_i,r_{i-1})$ from Rolle's Theorem and a root at all double roots of $g(x)$. The desired lemma follows easily.
\end{proof}
Now in order to prove that $P_{S_{n,t}}(x)$ has at most $2n-1$ roots which are roots of unity including multiplicity we instead study \[P_{S_{n,t}} (e^{in\theta})=e^{in\theta}(2\cos(n\theta)-2\cos((n-1)\theta)+\frac{\cos((2t+1)\theta)}{\cos(\theta)})\] and thus it suffices to bound the number of roots of \[Q_{n,t}(\theta)=(2\cos(n\theta)-2\cos((n-1)\theta)+\frac{\cos((2t+1)\theta)}{\cos(\theta)})\] in the interval $\left[0, 2\pi\right)$. 
\begin{lemma-hand}[5.]
$Q_{n,t}(\theta)$ has no roots in the intervals $\theta\in [0,\frac{1}{2t+4}]$ and $\theta\in \left[2\pi-\frac{1}{2t+4},2\pi\right)$.
\end{lemma-hand}
\begin{proof}
Due to symmetry around $\pi$, it suffices to prove the statement for the interval $\theta\in [0,\frac{1}{2t+4}]$. We consider two cases based on the parity of $t$. If $t=2k$ then by the sum-profuct formula for cosine,
\[\frac{\cos((2t+1)\theta)}{\cos(\theta)}=\sum_{j=1}^{k}2\cos(4j\theta)-2\cos((4j-2)\theta)+1.\]
Using this expansion it follows that if $\theta$ is a root that 
\[\Big|2\cos(n\theta)-2\cos((n-1)\theta)+1\Big|=\left|\sum_{j=1}^{k}2\cos(4j\theta)-2\cos((4j-2)\theta)\right|.\] Note that 
\begin{align*}
\Big|2\cos(n\theta)-2\cos((n-1)\theta)+1\Big|
&\ge 1-\Big|2\cos(n\theta)-2\cos((n-1)\theta)\Big|
\\&=1-4\left|\sin(\frac{\theta}{2})\sin(\frac{(2n-1)\theta}{2})\right|
\\&\ge1-4\Big|\sin(\frac{\theta}{2})\Big|.
\end{align*}
Furthermore note that 
\begin{align*}
\Big|\sum_{j=1}^{k}2\cos(4j\theta)-2\cos((4j-2)\theta)\Big|&\le 4\sum_{j=1}^{k}\Big|\sin\theta \sin((4j-1)\theta)\Big|
\\&\le 4k\Big|\sin\theta\Big|
\end{align*}
Combining these two estimates and the fact that $\Big|\sin(\theta)\Big|\le \Big|\theta\Big|$ with equality only at $0$ gives \[\Big|\theta\Big|> \frac{1}{4k+2}>\frac{1}{2t+4}.\] This finishes the case where $t$ is even. Otherwise $t=2k+1$ then 
\[\frac{\cos((2t+1)\theta)}{\cos(\theta)}=1+\sum_{j=0}^{k}2\cos((4j+2)\theta)-2\cos((4j)\theta).\]
Using this expansion it follows that if $\theta$ is a root that 
\[\Big|2\cos(n\theta)-2\cos((n-1)\theta)+1\Big|=\Big|\sum_{j=0}^{k}2\cos((4j+2)\theta)-2\cos(4j\theta)\Big|.\] As before note that 
\begin{align*}
\Big|2\cos(n\theta)-2\cos((n-1)\theta)+1\Big|
&\ge 1-\Big|2\cos(n\theta)-2\cos((n-1)\theta)\Big|
\\&=1-4\Big|\sin(\frac{\theta}{2})\sin(\frac{(2n-1)\theta}{2})\Big|
\\&\ge1-4\Big|\sin(\frac{\theta}{2})\Big|.
\end{align*}
Again note that 
\begin{align*}
\Big|\sum_{j=0}^{k}2\cos((4j+2)\theta)-2\cos(4j\theta)\Big|&\le 4\sum_{j=0}^{k}\Big|\sin\theta \sin((4j+1)\theta)\Big|
\\&\le 4(k+1)\Big|\sin\theta\Big|.
\end{align*}
Combining these two estimates and the fact that $\Big|\sin(\theta)\Big|\le \Big|\theta\Big|$ with equality only at $0$ gives \[\Big|\theta\Big|> \frac{1}{4k+6}=\frac{1}{2t+4}\] and therefore the desired estimate follows.
\end{proof}
\begin{lemma-hand}[6.]
For $n\ge \max\left(16(t+1)^3,40(t+2)\right)$, $Q_{n,t}(\theta)$ has at most $2n-1$ roots on the interval $\left[0,2\pi\right)$, including multiplicity.
\end{lemma-hand}
\begin{proof}
Now our strategy is to divide the remainder of the interval $\left[0, 2\pi\right)$ into at most $2n-2$ intervals of a fixed width, each of which contains exactly one root of $-\frac{Q_{n, t}'}{2}$.
For each positive integer $i\in \Big[(2n-1)\frac{1}{4\pi(t+2)}-1, (2n-1)(1-\frac{1}{4\pi(t+2)})\Big]$, we consider the intervals:
\[I_i=\left[\frac{(4i-1)\pi}{4n-2}, \frac{(4i+1)\pi}{4n-2}\right], J_i=\left[\frac{(4i+1)\pi}{4n-2}, \frac{(4i+3)\pi}{4n-2}\right],K_i= \left[\frac{(4i+3)\pi}{4n-2}, \frac{(4i+5)\pi}{4n-2}\right].\]
We claim that these intervals cover $[\frac{1}{2t+4}, 2\pi-\frac{1}{2t+4}]$. Indeed, since $I_{i+1}=K_i$, it suffices to check that these intervals cover some value at most $ \frac{1}{2t+4}$ and some value at least $ 2\pi-\frac{1}{2t+4}$. In particular note that an integer $i<(2n-1)\frac{1}{4\pi(t+2)}$ is included, and for this $i$ we have $\frac{(4i-1)\pi}{4n-2}<\frac{2i\pi}{2n-1}<\frac{1}{2t+4}$ as required. Similarly, there exists an included $i$ satisfying $i>(2n-1)(1-\frac{1}{4\pi(t+2)})-1$, so that $\frac{(4i+5)\pi}{4n-2}>\frac{2(i+1)\pi}{2n-1}>2\pi-\frac{1}{2t+4}$. Hence, all roots of $Q_n$ in $\left[0,2\pi\right)$ are in one of the $I_i, J_i$, or $K_i$.
\\
\\Next we claim that for each integer $i\in \Big[(2n-1)\frac{1}{4\pi(t+2)}-1, (2n-1)(1-\frac{1}{4\pi(t+2)})\Big]$, the following hold:
\begin{enumerate}[a.]
\item On the interval $I_i$, $\frac{-Q_{n, t}'}{2}$ and $(-1)^i$ have the same sign.
\item On the interval $J_i$, $\frac{-Q_{n, t}''}{2}$ is always nonzero.
\item On the interval $K_i$, $\frac{-Q_{n, t}'}{2}$ and $(-1)^i$ have opposite sign.
\end{enumerate}
Since $I_{i+1}=K_i$, we prove statements a and c simultaneously. It suffices to prove case a, including the extra interval $I_i$ for $i\in $($(2n-1)(1-\frac{1}{4\pi}), (2n-1)(1-\frac{1}{4\pi})+1$]. Note that we may rewrite $-\frac{Q_{n, t}'}{2}(\theta)$ as
\begin{align*}
-\frac{Q_{n, t}'}{2}(\theta)&= n\sin n\theta-(n-1)\sin ((n-1)\theta)+2\sum_{i=1}^t(-1)^{i+t}i\sin (2i\theta)
\\ &= 2n\sin\frac{\theta}{2}\cos \frac{(2n-1)\theta}{2}+\sin (n-1)\theta+2\sum_{i=1}^t(-1)^{i+t}i\sin (2i\theta)
\\ &=2n\sin\frac{\theta}{2}\cos \frac{(2n-1)\theta}{2}+ R_{n, t}(\theta)
\end{align*}
Note that the $R_{n, t}(\theta)$ terms sum to at most $1+2\sum_{i=1}^ti\le (t+1)^2$ in absolute value. On $I_i$, $\theta$ takes the form $\frac{(4i+\e)\pi}{4n-2}$ with $-1\le \e\le 1$. Now take $n\ge 40t+80$, so that $\frac{5\pi}{8n-4}< \frac{1}{16t+32}$. In particular
\[\frac{\theta}{2}\ge \frac{\pi}{8n-4}(4(\frac{2n-1}{4\pi(t+2)}-1)-1)=\frac{1}{4(t+2)}-\frac{5\pi}{8n-4}>\frac{3}{16(t+2)}.\]
Similarly
\[\frac{\theta}{2}\le \frac{\pi}{8n-4}(4((2n-1)(1-\frac{1}{4\pi(t+2)})+1)+1)=\pi-(\frac{1}{4(t+2)}-\frac{5\pi}{8n-4})<\pi-\frac{3}{16(t+2)}.\]
Hence $\sin\frac{\theta}{2}\ge \sin\frac{3}{16(t+2)}\ge \frac{1}{4\sqrt{2}(t+2)}$ for $\theta$ on all of these intervals. Furthermore, for $\theta=\frac{(4i+\e)\pi}{4n-2}$ we have
\[\cos\frac{(2n-1)\theta}{2}=\cos(i\pi+\e\frac{\pi}{4}).\]
Now we take two cases.
\\
\\ Case $1$: $i$ is even. Then $\cos(i\pi+\e\frac{\pi}{4})\ge \frac{\sqrt{2}}{2}$ and
\begin{align*}
\frac{-Q'_{n, t}}{2}&= 2n\sin\frac{\theta}{2}\cos \frac{(2n-1)\theta}{2}+R_{n, t}
\\ &> \frac{n}{4(t+2)}-(t+1)^2
\\ &> 0
\end{align*}
since $n\ge 8(t+1)^3$, as required.
\\
\\ Case $2$: $i$ is odd. Then $\cos(i\pi+\e\frac{\pi}{4})\le \frac{-\sqrt{2}}{2}$ and
\begin{align*}
\frac{-Q'_{n, t}}{2}&= -2n\sin\frac{\theta}{2}\cos \frac{(2n-1)\theta}{2}+R_{n, t}
\\ &< -\frac{n}{4(t+2)}+(t+1)^2
\\ &< 0
\end{align*}
since $n\ge 8(t+1)^3$, as required.
\\
\\ Now for case b we make use of our previous estimate on $\sin\frac{\theta}{2}$ for $\theta$ in $I_i$ and $K_i$, since the same estimate holds for the $J_i$. Note that
\begin{align*}-\frac{Q_{n, t}''(\theta)}{2}&=n^2\cos n\theta-(n-1)^2\cos (n-1)\theta+R_{n, t}'(\theta)=-2n^2\sin\frac{\theta}{2}\sin\frac{(2n-1)\theta}{2}+T_{n, t}(\theta).
\end{align*}
Here, $|T_{n, t}(\theta)|=|R'_{n, t}(\theta)+(2n-1)\cos (n-1)\theta|$, which is less than $3n-2+\sum_{i=1}^t 4i^2<3n-2+\frac{4}{3}(t+1)^3<4n$. 
On each $J_i$, the fractional part of $\frac{(2n-1)\theta}{2\pi}$ is between $\frac{1}{4}$ and $\frac{3}{4}$. Therefore it follows that $\Big|\sin\frac{(2n-1)\theta}{2}\Big|\ge \frac{\sqrt{2}}{2}$ and
\begin{align*}
\Big|2n^2\sin\frac{\theta}{2}\sin\frac{(2n-1)\theta}{2}\Big|&>\frac{n^2}{4(t+2)}
\\ &> 4n
\\ &>\Big|T_{n, t}(\theta)\Big|
\end{align*}
since $n\ge 40t+80$. So, $-\frac{Q_{n, t}''(\theta)}{2}$ is nonzero as required.
\\
\\ Now note that a, b, c together show that the only roots of $Q_n'(\theta)$ on $(\frac{1}{2t+4}, 2\pi-\frac{1}{2t+4})$ are contained in the $J_i$ and that none of these are roots of $Q_{n, t}''(\theta)$. This implies $Q_{n, t}(\theta)$ has no triple roots. Furthermore, for each $i$, we have that $Q_{n, t}'(\theta)$ takes opposite signs on $I_i, K_i$, and is monotonic on $J_i$. Hence, the number of roots of $Q'$ on this interval is at most the number of integer $i$ in $ [(2n-1)\frac{1}{4\pi(t+2)}-1, (2n-1)(1-\frac{1}{4\pi(t+2)})]$. Using the fact that $n\ge 40(t+2)$, this is at most
\begin{align*}(2n-1)(1-\frac{1}{4\pi(t+2)})-((2n-1)\frac{1}{4\pi(t+2)}-1)+1 &=2n-1-\frac{2n-1}{2\pi(t+2)}+2
\\ &<2n+1-\frac{n}{10(t+2)}
\\ &\le 2n+1-\frac{40(t+2)}{10(t+2)}
\\ &<2n-2.
\end{align*}

It follows that $Q_{n, t}'(\theta)$ has at most $2n-2$ roots in the interval with no double roots and therefore $Q_{n, t}(\theta)$ has at most $2n-1$ roots in $(\frac{1}{2t+4}, 2\pi-\frac{1}{2t+4})$ using Lemma 4. 
\end{proof}
Hence there are at most $2n-1$ roots of $P_{S_n}(x)$ which are roots of unity, so $S_n$ is not cyclotomic when $n\ge \max (40(t+2), 8(t+1)^3)$ as required for Theorem 1. The fact that the numerical semigroups $S_{n,t}$ corresponding to the ranges $t=0$, $5\le n\le 79$ and $t=1$, $8\le n\le 119$ are not cyclotomic may be verified directly via the numericalspgs GAP software \cite{gps}. See Appendix A. This verification gives all remaining dimensions $\ge 4$ for these constructions and therefore completes the proof of Theorem 2.

Since the polynomial $P_{S_{n,t}}(x)$ has a root off the unit circle, this naturally leads to the question of where this root off the unit circle is located, and the following result gives a bound on the variance of the roots from the unit circle. This result also appears to suggest the difficulty of demonstrating that there is a root off the unit circle for $P_{S_{n,t}}(x)$ given how close they are forced to lie to the unit circle for $t=0$.

\begin{thm-hand}[7.]
For positive $n\ge 12$, the roots of $P_{S_{n,0}}(x)=x^{2n}-x^{2n-1}+x^{n}-x+1$ satisfy $1-\frac{(\log n)^2}{n}\le |x|\le 1+\frac{(\log n)^2}{n}$.
\end{thm-hand}
\begin{proof}
Since $x$ being a root of $P_{S_{n,0}}(x)$ implies that $\frac{1}{x}$ is a root of $P_{S_{n,0}}(x)$ it suffices to prove the upper bound. Suppose for the sake of contradiction that $|x|> 1+\frac{(\log n)^2}{n}$. Then 
\begin{align*}|x|^{2n}&\ge \bigg(1+\frac{\log^2 n}{n}\bigg)|x|^{2n-1}
\\ &\ge |x|^{2n-1}+\frac{\log^2 n }{n}(1+\frac{\log^2 n}{n})^{n-1}|x|^{n}
\\ &\ge|x|^{2n-1}+\frac{\log^2 n}{n}((1+\frac{\log^2 n}{n})^{\frac{n}{\log^2 n }})^{\frac{\log^2 n}{2}}|x|^{n}
\\ &\ge|x|^{2n-1}+\frac{\log^2 n}{n}\Big(2^\frac{\log^2 n}{2}\Big)|x|^{n}
\\ &\ge |x|^{2n-1}+4|x|^n
\\ &>|x|^{2n-1}+|x|^2+|x|+1
\\ &\ge|x^{2n-1}-x^{n}+x-1|
\end{align*}
which is contradiction to the assumption that $x^{2n}-x^{2n-1}+x^{n}-x+1=0$ and the result follows. (Note that $n\ge 12$ is used third, fourth, and fifth inequalities.)
\end{proof}
This result, although surprising at first, becomes more natural in light of the result that for monic polynomials with coefficients chosen independently at random from $[-1,1]$, the roots lie near unit circle with high probability \cite{shepp1995complex}. Furthermore, Shepp and Vanderbei \cite{shepp1995complex} showed that in fact this occurs when the coefficients are chosen independently from the same distribution. 

\section{Conclusions}
In this paper, we've provided explicit constructions that demonstrate that the cyclotomic numerical semigroup condition is stronger than the symmetric numerical semigroup condition for all embedding dimensions at least $4$. In fact our results can be strengthened to show the density of roots in these polynomials which aren't roots of unity is asymptotically positive. It still remains open whether or not the cyclotomic property is equivalent to another natural condition, the \textit{complete intersection} property; this was also asked in \cite{Moree}.
\section{Acknowledgements}
The research was conducted at the University of Minnesota Duluth REU and was supported by NSF grant 1659047. The authors would like to thank Nathan Kaplan and Joe Gallian for suggesting the topic of numerical semigroups, Samuel Zbarsky for a significant simplification of the proof method, and Joe Gallian, Levent Alpoge, and Samuel Zbarsky for reading over the manuscript.
\bibliographystyle{plain}
\bibliography{main}

\begin{thebibliography}{1}

\bibitem{Moree}
Emil-Alexandru Ciolan, Pedro~A Garc{\'\i}a-S{\'a}nchez, and Pieter Moree.
\newblock Cyclotomic numerical semigroups.
\newblock {\em SIAM Journal on Discrete Mathematics}, 30(2):650--668, 2016.

\bibitem{Kronecker}
Pantellis~A. Damianou.
\newblock Monic polynomials in $\mathbb{Z}[x]$ with roots on the unit disc.
\newblock {\em American Math Monthly}, 108:253--257, 2001.

\bibitem{gps}
Manuel Delgado and Pedro~A Garc{\'\i}a-S{\'a}nchez.
\newblock numericalsgps, a gap package for numerical semigroups.
\newblock {\em ACM Communications in Computer Algebra}, 50(1):12--24, 2016.

\bibitem{moree2}
Herrera-Poyatos and Pieter Moree.
\newblock Higher order derivatives of the cyclotomic polynomial evaluated at
  $\pm 1$. (in preparation).

\bibitem{1}
Pieter Moree.
\newblock Numerical semigroups, cyclotomic polynomials, and bernoulli numbers.
\newblock {\em American Mathematical Monthly}, 121(10):890--902, 2014.

\bibitem{2}
Jos{\'e}~C Rosales and Pedro~A Garc{\'\i}a-S{\'a}nchez.
\newblock {\em Numerical semigroups}, volume~20.
\newblock Springer Science \& Business Media, 2009.

\bibitem{shepp1995complex}
Larry~A Shepp and Robert~J Vanderbei.
\newblock The complex zeros of random polynomials.
\newblock {\em Transactions of the American Mathematical Society}, pages
  4365--4384, 1995.

\end{thebibliography}
\section*{Appendix A}
The code for verifying the necessary cases for Theorem 2 is below. The numericalsgps package is obtained from \cite{gps}. 
\begin{lstlisting}[language=C++]
for n  in [5..79]  do
    Print(n, " ", IsCyclotomicNumericalSemigroup(
    NumericalSemigroup([n..(2*n-2)])), "\n");
    od;

\end{lstlisting}
\begin{lstlisting}[language=C++]
for n  in [8..119]  do
    Print(n, " ", IsCyclotomicNumericalSemigroup(
    NumericalSemigroup(Concatenation([n+2..(2*n-5)],[n-2, n-1]))), "\n");
    od;
\end{lstlisting}

\end{document}